\documentclass[11pt]{article}
\usepackage[a4paper,margin=25mm]{geometry}
\usepackage{amsmath,amssymb,amsthm}

\usepackage{xcolor}
\usepackage{hyperref}

\hypersetup{
    colorlinks=true,
    linkcolor=blue,   
    citecolor=red,    
    urlcolor=cyan!60!black   
}

\usepackage{setspace}
\newcommand{\R}{\mathbb R}
\newcommand{\Dc}[2]{\operatorname D^{\mathbf A}_{#1\to #2}}
\newcommand{\Ph}[2]{\Phi^{\mathbf A}_{#1\to #2}}
\newcommand{\Om}{\widetilde\Omega}
\newcommand{\Lc}{L^2_{\mathrm c}}
\newcommand{\esssupp}{\operatorname{ess\,supp}}
\numberwithin{equation}{section}
\title{Corrigendum to ``Well-posedness and regularity of the integral invariant model from linear scalar transport equation''\\[4pt]
\large [J. Math. Anal. Appl. 564 (2026) 130904]}
\author{Zhengrong Xie\\[3pt]
\small School of Mathematical Sciences, East China Normal University,\\
\small Shanghai 200241, China\\
\small \texttt{52265500018@stu.ecnu.edu.cn; xzr\_nature@163.com}}
\date{}
\begin{document}
\maketitle
\begin{abstract}
We clarify the compactly supported terminal test space and the moving-domain Bochner spaces in the integral invariant formulation. 
We complete the existence proof by establishing strong measurability and the initial trace, and replace the stability argument by a compact-test density argument that retains the original stability constant. 
A short appendix shows that the integral invariant identity remains valid for all square-integrable terminal test functions that vanish almost everywhere outside the terminal domain, thereby admitting the zero extensions of the cell polynomials used in the SLDG/SLFV method.  
We also correct the velocity and divergence in the swirling-deformation example. 

\end{abstract}

\tableofcontents
\clearpage

All references to original page, equation and theorem numbers below concern \cite{original}. We retain the notation
\[
H(\R^d)=L^2(\R^d),\qquad
\Om(t)=\Dc{T}{t}(\Omega),\qquad
\Om^0=\Om(0),
\]
where $\Omega\subset\R^d$ is bounded and open.  
The associated space-time domain is denoted by
\[
\widetilde Q_T
:=
\{(\mathbf x,t):t\in[0,T],\ \mathbf x\in\Om(t)\}.
\]
The flow maps are denoted by $\Dc{s}{t}$, their spatial Jacobian matrices by $J_{\Dc{s}{t}}$, and the test-function evolution by
\[
\Ph{s}{t}[f](\mathbf x)=f\bigl(\Dc{t}{s}(\mathbf x)\bigr).
\]
All scalar functions and Hilbert spaces in this notice are real-valued.

\section{Function spaces and the weak formulation}
\label{sec:spaces}

\paragraph{Compactly supported terminal tests.}
The space denoted by $H_0(\Omega)$ on p.~5 is generally not complete in the $L^2$ norm and should not be called a Hilbert space. We write it as
\begin{equation}\label{eq:compact}
\Lc(\Omega):=
\{f\in L^2(\R^d):\esssupp f\subset\subset\Omega\}.
\end{equation}
Here $\subset\subset$ means compact containment, and support is understood as essential support for $L^2$ equivalence classes. This is an inner-product subspace of $H(\R^d)$. We retain $\Lc(\Omega)$ as the terminal test space in the corrected definition and proofs below. Its completeness is not needed: the Riesz representation theorem is applied in the full Hilbert space $H(\R^d)$. 
Appendix~\ref{app:extended} shows, by density, that the integral
invariant identity remains valid when the terminal test space is enlarged from $\Lc(\Omega)$ to $\mathcal H_\Omega$, 
the space of $L^2(\R^d)$ functions vanishing almost everywhere outside $\Omega$. 
When $\Omega$ is a mesh cell, this larger space contains the zero extensions of cell polynomials and thus accommodates
the terminal tests used in SLDG/SLFV numerical methods. 

\paragraph{Moving-domain Bochner spaces.}
For an open set $D$, let $E_D:L^2(D)\to L^2(\R^d)$ be zero extension, and put $E_t=E_{\Om(t)}$. This map is an isometry. In Definition 2.1, strong measurability must be required in addition to the norm condition. For the moving domains in Definition 3.1 and Remark 3.1, we use
\begin{equation}\label{eq:bochner}
U\in L^q([0,T];H(\Om(t)))
\quad\Longleftrightarrow\quad
\bigl[t\mapsto E_tU(\cdot,t)\bigr]
\in L^q([0,T];H(\R^d)),\qquad 1\leq q\leq\infty.
\end{equation}
The spatial slices are $L^2$ elements for almost every time, and families are identified when their zero extensions agree almost everywhere in time. Strong measurability means almost-everywhere approximation in the Hilbert-space norm by measurable simple functions. The norms are
\begin{align}
\|U\|_{L^q([0,T];H(\Om(t)))}
&=\left(\int_0^T\|U(\cdot,t)\|_{L^2(\Om(t))}^q\,dt\right)^{1/q},
&&1\leq q<\infty,\label{eq:qnorm}\\
\|U\|_{L^\infty([0,T];H(\Om(t)))}
&=\operatorname*{ess\,sup}_{0<t<T}\|U(\cdot,t)\|_{L^2(\Om(t))}.
\label{eq:infnorm}
\end{align}
In Remark~2.3, the case $q=\infty$ is understood in the
standard Bochner sense, requiring strong measurability and
essential boundedness of the spatial norm.
Lemma~4.4, applied to this scalar norm, provides the
characterization through uniformly bounded finite-$q$
norms and the corresponding norm-limit identity. 

\paragraph{Flow estimates.}
We keep the coefficient assumptions in Section 2 of the original article. The $C^1$ velocity field has a globally bounded spatial derivative and satisfies $|\nabla\cdot\mathbf A|\leq M_A$. The spatial Lipschitz bound and boundedness of $\mathbf A(0,t)$ on $[0,T]$ give linear growth and hence global $C^1$ flow diffeomorphisms. Liouville's formula gives
\begin{align}
\det J_{\Dc{s}{t}}(\mathbf x)
&=\exp\left(\int_s^t\nabla\cdot\mathbf A
\bigl(\Dc{s}{\tau}(\mathbf x),\tau\bigr)\,d\tau\right)>0,
\label{eq:jac}\\
e^{-M_A|t-s|}&\leq\det J_{\Dc{s}{t}}(\mathbf x)
\leq e^{M_A|t-s|}.\label{eq:jacbound}
\end{align}
The flow and its inverse preserve null sets, so composition
with these maps is independent of the chosen representatives
of the $L^2$ classes.
For $s,t\in[0,T]$ and $f\in\Lc(\Om(s))$,
the support-transport argument of Lemma~4.1,
interpreted in terms of essential support and applied
from $s$ to $t$, gives
\begin{equation}\label{eq:support}
\esssupp\bigl(\Ph{s}{t}[f]\bigr)
=
\Dc{s}{t}(\esssupp f)
\subset\subset\Om(t).
\end{equation}
Applying the norm estimate of Lemma~4.2 with $P=2$
to the same evolution yields
\begin{equation}\label{eq:flowbound}
\|\Ph{s}{t}[f]\|_{L^2(\Om(t))}
\leq
e^{M_A|t-s|/2}\|f\|_{L^2(\Om(s))}.
\end{equation}
Together, \eqref{eq:support} and \eqref{eq:flowbound} imply
\[
\Ph{s}{t}[f]\in\Lc(\Om(t)).
\]
The change-of-variables proof of \eqref{eq:flowbound} applies to every $f\in L^2(\Om(s))$; for whole-space compositions, its zero extension is understood.

\paragraph{Precise weak formulation.}
For $\widetilde U_0\in L^2(\Om^0)$, Definition 3.1 is understood as follows. A family $U\in L^1([0,T];H(\Om(t)))$ is an integral-invariant weak solution if it admits a representative with $U(\cdot,t)\in L^2(\Om(t))$ for every $t\in[0,T]$, such that
\begin{equation}\label{eq:invariant}
\int_{\Om(t)}U(\mathbf x,t)\Ph{T}{t}[\Psi](\mathbf x)\,d\mathbf x
=\int_{\Om^0}\widetilde U_0(\mathbf x)\Ph{T}{0}[\Psi](\mathbf x)\,d\mathbf x
\end{equation}
holds for every $t\in[0,T]$ and every $\Psi\in\Lc(\Omega)$ on this representative. This is the precise integrated form of (3.3b).

Equivalently, in the differential formulation, require the scalar pairings to be absolutely continuous on $[0,T]$, to have derivative zero almost everywhere, and to have the initial values prescribed by the right-hand side of \eqref{eq:invariant}. The equivalence follows from the fundamental theorem for absolutely continuous functions. The initial condition is therefore a condition on the initial values of these pairings, not an arbitrary assignment at a single time of a Bochner equivalence class. Pulling back arbitrary $C^\infty_{\mathrm c}(\Om^0)$ functions and using density gives $U(\cdot,0)=\widetilde U_0$ in $L^2(\Om^0)$. The strong initial trace will follow from the existence construction and uniqueness.

\section{Completion of the existence proof}
\label{sec:existence}

We complete the proof of Theorem 5.1 while retaining
the functional (5.6) and its whole-space Riesz construction,
with $t^*$ denoted by $\lambda$.
For ease of reference, we recall the Riesz identity (5.7)
and the bound for $\mathcal W^\lambda$ contained in (5.20):
\begin{align}
\langle\mathcal W^\lambda,w\rangle_{L^2(\R^d)}
&=
f^0(w;\widetilde U_0,\mathbf A,0,\lambda),
\qquad w\in H(\R^d),
\label{eq:riesz}\\
\|\mathcal W^\lambda\|_{L^2(\R^d)}
&\leq
e^{M_A\lambda/2}
\|\widetilde U_0\|_{L^2(\Om^0)}.
\label{eq:bound}
\end{align}
The additional arguments below identify $\mathcal W^\lambda$
with the zero extension of the solution and establish
strong measurability in time and the strong initial trace.

For general $\widetilde U_0\in L^2(\Om^0)$,
only its zero extension to $\R^d$ is required.
The compact-support assertions after (5.1) and in (5.19)
should therefore be replaced by this zero-extension property.

Fix $\lambda\in[0,T]$ and define
$U(\cdot,\lambda)=\mathcal W^\lambda|_{\Om(\lambda)}$.
Since \eqref{eq:riesz} holds for every $w\in H(\R^d)$,
we may take
\[
w=\mathcal W^\lambda
\chi_{\R^d\setminus\Om(\lambda)}.
\]
Its right-hand side vanishes because
$\Dc{0}{\lambda}(\Om^0)=\Om(\lambda)$, and hence
\[
\int_{\R^d\setminus\Om(\lambda)}
|\mathcal W^\lambda|^2\,d\mathbf x=0.
\]
Consequently,
\begin{equation}\label{eq:extension}
E_\lambda U(\cdot,\lambda)=\mathcal W^\lambda,
\qquad
\|U(\cdot,\lambda)\|_{L^2(\Om(\lambda))}
=
\|\mathcal W^\lambda\|_{L^2(\R^d)}.
\end{equation}
This clarifies the relation between the restriction in (5.14)
and its zero extension.
The restricted tests in (5.10), (5.13) and (5.22)
determine only the restriction of a whole-space function
to the corresponding spatial domain.
Uniqueness of $\mathcal W^\lambda$ on the whole space follows
from \eqref{eq:riesz}.

For any $\Psi\in\Lc(\Omega)$, set
$\psi(\cdot,\lambda)=\Ph{T}{\lambda}[\Psi]$.
By \eqref{eq:support}--\eqref{eq:flowbound},
this is an admissible choice of $w$ in \eqref{eq:riesz}.
The flow composition law gives
\[
\psi\bigl(\Dc{0}{\lambda}(\mathbf y),\lambda\bigr)
=
\psi(\mathbf y,0)
\quad\text{for a.e. }\mathbf y\in\Om^0.
\]
Equations \eqref{eq:riesz} and \eqref{eq:extension}
therefore give \eqref{eq:invariant}
for every compact terminal test and every time.

\paragraph{Strong measurability.}
For $\varphi\in C^\infty_{\mathrm c}(\R^d)$,
equation \eqref{eq:riesz} shows that
$\lambda\mapsto\langle\mathcal W^\lambda,\varphi\rangle$
is continuous.
Indeed, the flow is continuous in time,
$\widetilde U_0\in L^1(\Om^0)$ because $\Om^0$ is bounded,
and Lebesgue dominated convergence applies with majorant
$|\widetilde U_0|\|\varphi\|_\infty$.
For an arbitrary $w\in L^2(\R^d)$,
choose $\varphi_n\in C^\infty_{\mathrm c}(\R^d)$
converging to $w$ in $L^2$.
The uniform bound
\[
\sup_{0\leq\lambda\leq T}
|\langle\mathcal W^\lambda,w-\varphi_n\rangle|
\leq
e^{M_AT/2}\|\widetilde U_0\|_{L^2(\Om^0)}
\|w-\varphi_n\|_{L^2(\R^d)}
\longrightarrow0
\]
proves \textbf{weak continuity of
$\lambda\mapsto\mathcal W^\lambda$.
Since $L^2(\R^d)$ is separable,
the Pettis measurability theorem gives strong measurability}.
Together with \eqref{eq:bound}--\eqref{eq:extension},
this fills the gap between (5.20) and (5.21)
and proves that the constructed family belongs
to the required Bochner space.

\paragraph{Initial trace and temporal integrability.}
At $\lambda=0$, \eqref{eq:riesz} gives
$\mathcal W^0=E_0\widetilde U_0$.
Writing $\|\cdot\|_2$ for the norm in $L^2(\R^d)$,
the isometry of zero extension and \eqref{eq:bound} yield
\(
\|\mathcal W^\lambda\|_2^2
\leq e^{M_A\lambda}\|\mathcal W^0\|_2^2.
\)
Moreover, weak continuity gives
\(
\langle\mathcal W^\lambda,\mathcal W^0\rangle_{L^2(\R^d)}
\longrightarrow
\|\mathcal W^0\|_2^2
\quad\text{as }\lambda\downarrow0.
\)
Consequently,
\begin{align*}
0
\leq \|\mathcal W^\lambda-\mathcal W^0\|_2^2
&=
\|\mathcal W^\lambda\|_2^2+\|\mathcal W^0\|_2^2
-2\langle\mathcal W^\lambda,\mathcal W^0\rangle_{L^2(\R^d)}\\
&\leq
\bigl(e^{M_A\lambda}+1\bigr)\|\mathcal W^0\|_2^2
-2\langle\mathcal W^\lambda,\mathcal W^0\rangle_{L^2(\R^d)}
\longrightarrow0
\qquad\text{as }\lambda\downarrow0.
\end{align*}
By the squeeze theorem and \eqref{eq:extension},
the constructed solution has the strong initial trace
\begin{equation}\label{eq:trace}
\lim_{\lambda\downarrow0}
\|E_\lambda U(\cdot,\lambda)-E_0\widetilde U_0\|_{L^2(\R^d)}
=0.
\end{equation}
This replaces the justification by a single-time redefinition
in Step 4 of the original proof. 

The same measurability argument and norm bound
also give Theorem 6.1 directly:
\begin{align}
\|U\|_{L^\infty([0,T];H(\Om(\lambda)))}
&\leq
e^{M_AT/2}\|\widetilde U_0\|_{L^2(\Om^0)},
\label{eq:linfty}\\
\|U\|_{L^q([0,T];H(\Om(\lambda)))}
&\leq
T^{1/q}e^{M_AT/2}\|\widetilde U_0\|_{L^2(\Om^0)},
\qquad 1\leq q<\infty.
\label{eq:lq}
\end{align}
The norm-limit assertion (6.2) follows from the corresponding
finite-measure identity for the bounded measurable scalar
function
$\lambda\mapsto\|U(\cdot,\lambda)\|_{L^2(\Om(\lambda))}$.

\section{Stability with compact terminal tests and uniqueness}
\label{sec:stability}

The argument in (5.32)--(5.38) requires correction. The terminal test formed from the solution difference in (5.33) was not shown to belong to the compact test class, and the change of variables in (5.35) omitted a Jacobian. 
We replace that argument by the following proof using only admissible compact tests. 
Appendix~\ref{app:extended} explains when the original test constructed from the difference of two weak solutions can also be used. 

Let $U,V$ be any two weak solutions in the sense of Section~\ref{sec:spaces}, for the same velocity field and with respective initial data $\widetilde U_0,\widetilde V_0\in L^2(\Om^0)$. Use their representatives satisfying \eqref{eq:invariant} at every time and write
\[
w_t=U(\cdot,t)-V(\cdot,t),\qquad
\widetilde w_0=\widetilde U_0-\widetilde V_0,\qquad
J_t(\mathbf y)=\det J_{\Dc{0}{t}}(\mathbf y)>0.
\]
Fix $t\in[0,T]$. For any $\varphi\in C^\infty_{\mathrm c}(\Om^0)$, identified here with its smooth zero extension to $\R^d$, choose
\begin{equation}\label{eq:compacttest}
\Psi(\mathbf z)=\varphi\bigl(\Dc{T}{0}(\mathbf z)\bigr),
\qquad \mathbf z\in\R^d.
\end{equation}
Equations \eqref{eq:support}--\eqref{eq:flowbound} ensure $\Psi\in\Lc(\Omega)$. Its evolution satisfies
\[
\psi(\mathbf x,\lambda)
=\varphi\bigl(\Dc{\lambda}{0}(\mathbf x)\bigr),
\qquad \mathbf x\in\R^d,\quad \lambda\in[0,T],
\]
and in particular $\psi(\cdot,0)=\varphi$. Subtracting \eqref{eq:invariant} for $U,V$ and changing variables $\mathbf x=\Dc{0}{t}(\mathbf y)$ gives
\begin{equation}\label{eq:densitypair}
\int_{\Om^0}
\left[w_t\bigl(\Dc{0}{t}(\mathbf y)\bigr)J_t(\mathbf y)
-\widetilde w_0(\mathbf y)\right]\varphi(\mathbf y)\,d\mathbf y=0.
\end{equation}
The expression in brackets belongs to $L^2(\Om^0)$, since \eqref{eq:jacbound} yields
\begin{align*}
\int_{\Om^0}|w_t(\Dc{0}{t}(\mathbf y))J_t(\mathbf y)|^2\,d\mathbf y
\leq e^{M_At}\int_{\Om^0}|w_t(\Dc{0}{t}(\mathbf y))|^2J_t(\mathbf y)\,d\mathbf y
=e^{M_At}\|w_t\|_{L^2(\Om(t))}^2.
\end{align*}
By the density of $C^\infty_{\mathrm c}(\Om^0)$ in $L^2(\Om^0)$, \eqref{eq:densitypair} implies
\begin{equation}\label{eq:densityidentity}
w_t\bigl(\Dc{0}{t}(\mathbf y)\bigr)J_t(\mathbf y)
=\widetilde w_0(\mathbf y)
\quad\text{for a.e. }\mathbf y\in\Om^0.
\end{equation}
Another application of the change-of-variables formula gives
\[
\|w_t\|_{L^2(\Om(t))}^2
=\int_{\Om^0}\frac{|\widetilde w_0(\mathbf y)|^2}{J_t(\mathbf y)}\,d\mathbf y
\leq e^{M_At}\|\widetilde w_0\|_{L^2(\Om^0)}^2.
\]
Thus
\begin{equation}\label{eq:stability}
\|U(\cdot,t)-V(\cdot,t)\|_{L^2(\Om(t))}
\leq e^{M_At/2}\|\widetilde U_0-\widetilde V_0\|_{L^2(\Om^0)},
\qquad 0\leq t\leq T.
\end{equation}
In particular, (5.27) and (5.38) hold with the original uniform constant $e^{M_AT/2}$. The proof applies to arbitrary weak solutions and arbitrary $L^2$ initial data; the compact-support restriction on the initial data in Theorem 5.2 is unnecessary.

For equal initial data, \eqref{eq:stability} gives equality in $L^2(\Om(t))$ at every time. This proves Theorem 5.3 and equality of the corresponding moving-domain Bochner elements. To interpret the conclusion on $\widetilde Q_T$, choose jointly measurable representatives of the zero extensions. Their slices represent the original $L^2$ elements for almost every time, and Tonelli's theorem gives
\[
\int_0^T\int_{\R^d}|E_tU(\mathbf x,t)-E_tV(\mathbf x,t)|^2\,d\mathbf x\,dt=0.
\]
Thus $U=V$ almost everywhere on $\widetilde Q_T$. Every weak solution consequently coincides with the constructed solution and inherits \eqref{eq:trace}--\eqref{eq:lq}.

\section{Correction of the swirling-deformation example}

In the third two-dimensional example in Appendix B (p.~20), the displayed PDE has the intended sign in its first flux, but the first component of the displayed velocity is missing a minus sign. The velocity and its spatial Jacobian should read
\begin{align}
\mathbf A(x,y,t)
&=\cos\left(\frac{\pi t}{T}\right)
\begin{pmatrix}
-2\pi\cos^2(x/2)\sin y\\
2\pi\sin x\cos^2(y/2)
\end{pmatrix},\label{eq:swirl}\\
J_{\mathbf A}(x,y,t)
&=\cos\left(\frac{\pi t}{T}\right)
\begin{pmatrix}
\pi\sin x\sin y&-2\pi\cos^2(x/2)\cos y\\
2\pi\cos^2(y/2)\cos x&-\pi\sin x\sin y
\end{pmatrix}.
\label{eq:swirljac}
\end{align}
In particular,
\begin{equation}\label{eq:swirldiv}
\nabla\cdot\mathbf A=0,\qquad
\sup_{(x,y,t)\in\R^2\times[0,\infty)}|\nabla\cdot\mathbf A|=0.
\end{equation}
The printed divergence supremum $2\pi$ must therefore be replaced by $0$.

\paragraph{Effect of the corrections.}
With the original coefficient assumptions and the weak formulation in Section~\ref{sec:spaces}, the corrected arguments establish existence, stability, uniqueness and temporal integrability. The general stability constant is unchanged. The swirling-deformation example has the corrected velocity \eqref{eq:swirl} and zero divergence. The appendix derives the extension of the integral identity needed for local DG tests.

\appendix
\section{Extension of the test class and cellwise DG tests}
\label{app:extended}

The integral invariant model originates from the cellwise construction of the SLDG/SLFV method: on each time step, a terminal mesh cell is traced back along the characteristics, and a polynomial on that terminal cell supplies the terminal test. For a bounded open mesh cell $K$ and a polynomial $p\in\mathbb P^k(K)$, the zero extension $E_Kp$ belongs to $L^2(\R^d)$ but need not belong to $\Lc(K)$. For example, $E_K1=\chi_K$ has essential support $\overline K$. This motivates the larger space
\begin{equation}\label{eq:hspace}
\mathcal H_\Omega
:=\{f\in L^2(\R^d):f=0\text{ a.e. on }\R^d\setminus\Omega\}
=E_\Omega L^2(\Omega).
\end{equation}
Zero extension is an isometric isomorphism from $L^2(\Omega)$ onto $\mathcal H_\Omega$, so the latter is a closed Hilbert subspace of $H(\R^d)$. Its elements need not have essential support compactly contained in $\Omega$, and no Sobolev boundary trace is imposed.

\paragraph{Extension of the integral identity.}
Let $U$ be any weak solution from Section~\ref{sec:spaces}, with the representative already chosen there. For $\Psi\in\mathcal H_\Omega$, choose zero-extended $\Psi_m\in C^\infty_{\mathrm c}(\Omega)$ such that $\Psi_m\to\Psi$ in $L^2(\R^d)$. For each fixed $t$, estimate \eqref{eq:flowbound} gives
\begin{align*}
&\left|\int_{\Om(t)}U(\mathbf x,t)\Ph{T}{t}[\Psi-\Psi_m](\mathbf x)\,d\mathbf x\right|
+\left|\int_{\Om^0}\widetilde U_0(\mathbf x)\Ph{T}{0}[\Psi-\Psi_m](\mathbf x)\,d\mathbf x\right|\\
&\leq\left(e^{M_A(T-t)/2}\|U(\cdot,t)\|_{L^2(\Om(t))}
+e^{M_AT/2}\|\widetilde U_0\|_{L^2(\Om^0)}\right)
\|\Psi-\Psi_m\|_{L^2(\R^d)}\longrightarrow0.
\end{align*}
Passing to the limit in \eqref{eq:invariant} for $\Psi_m$ proves that the same identity holds for every $\Psi\in\mathcal H_\Omega$ and every $t\in[0,T]$, on the same representative. Its scalar pairings are again constant and hence absolutely continuous. Conversely, requiring all $\mathcal H_\Omega$ tests includes the compact tests. Therefore this extension requires no second existence or uniqueness proof.

Applying the same formulation to a time step $[t^n,t^{n+1}]$ with terminal domain $K$, every $p\in\mathbb P^k(K)$ is admitted through $E_Kp$, giving the continuous cell identity
\begin{equation}\label{eq:cell}
\int_K U(\mathbf x,t^{n+1})p(\mathbf x)\,d\mathbf x
=\int_{\Dc{t^{n+1}}{t^n}(K)}
U(\mathbf y,t^n)
p\bigl(\Dc{t^n}{t^{n+1}}(\mathbf y)\bigr)\,d\mathbf y.
\end{equation}
This supplies the terminal tests used in the local integral relation underlying SLDG/SLFV. It concerns the continuous model; numerical characteristic tracing, quadrature and projection remain separate discretization steps.

\paragraph{The solution-difference test.}
The extension also justifies a corrected version of the original stability argument. Fix $t\in[0,T]$ and set
\begin{equation}\label{eq:difftest}
\Psi_t(\mathbf z)
:=\bigl(E_t[U(\cdot,t)-V(\cdot,t)]\bigr)
\bigl(\Dc{T}{t}(\mathbf z)\bigr),
\qquad \mathbf z\in\R^d.
\end{equation}
Since $\Dc{T}{t}$ maps $\R^d\setminus\Omega$ onto
$\R^d\setminus\Om(t)$, and
$E_t[U(\cdot,t)-V(\cdot,t)]$ vanishes almost everywhere
outside $\Om(t)$, preservation of null sets under the flow gives
\[
\Psi_t=0
\quad\text{a.e. on }\R^d\setminus\Omega.
\]
Applying \eqref{eq:flowbound} to $U(\cdot,t)-V(\cdot,t)$
for the evolution from $t$ to $T$ also gives
$\Psi_t\in L^2(\R^d)$.
Hence $\Psi_t\in\mathcal H_\Omega$. 
The corrected change of variables in (5.35) is
\begin{equation}\label{eq:missingjac}
\|\Psi_t\|_{L^2(\Omega)}^2
=\int_{\Om(t)}|U(\mathbf x,t)-V(\mathbf x,t)|^2
\det J_{\Dc{t}{T}}(\mathbf x)\,d\mathbf x.
\end{equation}
Once $t$ is fixed, use $\lambda$ for the evolution time. The flow composition law yields
\begin{equation}\label{eq:diffevolution}
\psi_t(\cdot,\lambda):=\Ph{T}{\lambda}[\Psi_t]
=\Ph{t}{\lambda}\bigl[E_t[U(\cdot,t)-V(\cdot,t)]\bigr].
\end{equation}
In particular, its restriction at $\lambda=t$ equals $U(\cdot,t)-V(\cdot,t)$ in $L^2(\Om(t))$. Applying \eqref{eq:flowbound} directly from $t$ to $0$ gives
\[
\|\psi_t(\cdot,0)\|_{L^2(\Om^0)}
\leq e^{M_At/2}\|U(\cdot,t)-V(\cdot,t)\|_{L^2(\Om(t))}.
\]
Subtracting the now extended invariant identities at time $t$ and applying Cauchy--Schwarz therefore gives
\begin{align*}
\|U(\cdot,t)-V(\cdot,t)\|_{L^2(\Om(t))}^2
&=\int_{\Om^0}(\widetilde U_0-\widetilde V_0)(\mathbf y)
\psi_t(\mathbf y,0)\,d\mathbf y\\
&\leq e^{M_At/2}\|\widetilde U_0-\widetilde V_0\|_{L^2(\Om^0)}
\|U(\cdot,t)-V(\cdot,t)\|_{L^2(\Om(t))}.
\end{align*}
This recovers \eqref{eq:stability}. The test in \eqref{eq:difftest} is used only after its admissibility has been established through \eqref{eq:hspace} and the density extension above.

\end{document}